# The *q*-Gauss-Newton method for unconstrained nonlinear optimization

Danijela D. Protić
University of Nis
Center for Applied Mathematics and Electronic
Belgrade 11090
adanijela@ptt.rs

Miomir M. Stanković
Mathematical Institute of SASA
Belgrade 11000
miomir.stankovic@gmail.com

## ABSTRACT

A *q*-Gauss-Newton algorithm is an iterative procedure that solves nonlinear unconstrained optimization problems based on minimization of the sum squared errors of the objective function residuals. Main advantage of the algorithm is that it approximates matrix of *q*-second order derivatives with the first-order *q*-Jacobian matrix. For that reason, the algorithm is much faster than *q*-steepest descent algorithms. The convergence of *q*-GN method is assured only when the initial guess is close enough to the solution. In this paper the influence of the parameter *q* to the non-linear problem solving is presented through three examples. The results show that the *q*-GD algorithm finds an optimal solution and speeds up the iterative procedure.

## 1 Introduction

In his work Theoria motus, Carl Friedrich Gauss predicted position of reappearing of the dwarf planet Ceres, first slighted in 1801, accurate within half a degree [1]. In the following years Gauss developed a more precise method of calculating the effect of the Sun and larger planets on the orbits of planetoids by reducing inaccuracy of calculations arising from measurement errors. The primacy of this discovery was disputed between him and Legendre who's Essay of the numbers' theory was also reissued in this series [2]. Gauss' iterative algorithm is used to solve systems of linear equations and find inverse an invertible matrix. The fundamental idea is to add multipliers of one equation to the others in order to eliminate a variable and to continue this process until only one variable is left. Ones the final variable is determined, it is substituted backwards into the other equations in order to eliminate the remaining unknowns [3]. In decades many methods for solving the nonlinear optimization problems based on Gauss' theory are developed. The most known are steepest-descent (SD) [4], [5], [6], quasi-Newton (QN) [7], [8] and conjugate-gradient (CG) methods [9], [10], [11]. The main idea of the SD method is to minimize the objective function by the iterative procedure based on the step length and the search direction that is determined by the negative of the gradient. The SD uses the first order derivative of objective function to find local and global minima. The training process of the SD algorithm is asymptotic convergence. Around the solution all the elements of gradient vector are very small and there is very tiny function change. QN methods require only the gradient of the objective function to be supplied at each iterate and are often used to solve a system of nonlinear equations. CG minimization is recognized as one of the best methods in the field of unconstrained optimization due to its simplicity and low memory management [12]. The Newton's method is much faster than all of the abovementioned methods but has a quadratic order of convergence and requires a lot of computational works. Comparing with the SD method, the Hessian matrix of second-order derivatives of the objective function need to be calculated for each component of the gradient vector. The Gauss-Newton (GN) method is an approximation of the Newton's method [13]. Main advantage of the GN algorithm over the standard Newton's method is that former does not require calculation of second-order derivatives by introducing the Jacobian matrix in order to simplify the process of calculating the Hessian matrix directly. Thus, only the first-order terms need to be evaluated [14], [15], [16]. The GN method is widely used for various purposes such as inversion of acoustic data [17], regularization in microwave biomedical imaging [18], cancer classification [19], optimization of the hydrological models [20] and reduction of the negative effects of numerical noise and truncation errors [21]. As the algorithm presumes that the error function is approximately quadratic near the optimal solution, it is frequently used for solving the weighted non-linear least squares (WNLS) problems [22] and approximating the sum squared errors (SSE) in neural network optimization [23], [24], [25], [26]. However, the GN algorithm still faces same convergent problems like the Newton algorithm for complex optimization [27].

Over the years various methods for solving unconstrained optimization problems based on $q$-calculus have been proposed. The $q$-calculus is widely used for solving problems in the fields of number theory [28], [29], physics [30], [31], quantum mechanics [32], engineering [33] and economics [34], [35]. Euler and Hine obtained basic formulae in $q$-calculus [12] while Frank Hilton Jackson introduced $q$-derivative and $q$-integral [36]. Classical SD method and $q$-SD method are applied to unimodal and multimodal test functions in [37]. The results have shown great performance of the $q$-SD method. The algorithm was able to escape from many local minima and reach the global minimum for multimodal function. The $q$-SD developed to solve objective unconstrained optimization problems where the $q$-parameter acts as the dilation parameter is given in [21]. The proposed algorithm is complemented with strategies for generating the parameter $q$ and computing the step length. Given strategies are implemented so that the search process gradually shifts from global minimum at the beginning to local minimum as the algorithm converges and approaches the optimal value [38]. In [33] the $q$-versions of the QN models ($q$-QN) are presented. Authors adapted Newton's, Newton-Kantorowich and gradient methods to cases when the functions are given in the form of infinitive products. Given examples comprehended the infinitive $q$-power products and proved that the proposed methods were suitable for solving the problems. In [39] a limited memory Broyden-Fletcher-Glodfarb-Shanno (BFGS) method for solving optimization problems using $q$-derivatives is described. The Jackson's derivative is used as a mechanism for escaping from local minima, while the $q$-gradient method is complemented to generate the parameter $q$ for computing the step length. The $q$-CG model given in [40] proved the convergence of the $q$-CG method without assumption of convexity. In [41] the $q$-CG model is used with the $q$-Fletcher-Reeves ($q$-FR) method for solving unconstrained optimization problems. The $q$-descent direction of the objective function is derived using $q$-gradient, a generalization of the classical gradient based on the $q$-derivative. Authors proved that the proposed algorithm was globally convergent with the Armijo-type line search condition. The algorithm solves the equation by using damped Newton iteration with the line search strategy to improve convergence from bad initial guess.

In this paper we present a $q$-Gauss-Newton ($q$-GN) algorithm based on implemented $q$-first-order derivatives of the components of the objective function. The $q$-GN algorithm is iterative, local and converges only when the iteration process starts near to the optimal solution. The algorithm does not require evaluations of the $q$-second-order derivatives since the $q$-Hessian matrix is approximated with the $q$-Jacobian matrix of the $q$-first-order derivatives of the objective function. An SSE as well as the number of iterations are used as stopping criteria. The influence of the parameter $q$ to the nonlinear problem solving is presented through three examples. First example demonstrates finding a global minimum of the nonlinear function. The results based on different values of the parameter $q$ are compared with the result based on Nelder-Mead simplex algorithm which finds minima by using damped Newton's method [42], [43]. The maximum number of the Gauss-Newton iteration steps is set to 25. The initial guess is chosen to be near global minimum. In the second example the optimal value of the real-valued function of two variables $f: R^2 \rightarrow R^2$ is determined. It is shown that the iterative procedure slightly changes the results with increasing the parameter $q$. The optimal solution of a real-valued function of two variables $f: R^2 \rightarrow R^3$ is determined in the third example. The results show that number of iterations decrease gradually as the parameter $q$ increases. The examples show that the $q$-GD algorithm finds a global minimum of the objective function and speeds up the iterative procedure near the optimal solution.

## 2  $q$-calculus preliminarities

Let the $q$-analog of $n \in N$, be defined as:

$$[n] = \frac{q^n - 1}{q - 1} \tag{1}$$

where $0 < q < 1$. Then the $q$-binomial coefficients and $q$-factorial of a $[n]$ can be expressed as follows [44]:

$$\begin{bmatrix} n \\ j \end{bmatrix} = \frac{[n]!}{[j]![n-j]!} \tag{2}$$

### 2.1  $q$-differential and $q$-derivative

Consider an arbitrary function $f(x)$. The $q$-differential of $f(x)$ is given with [44]:

$$d_q f(x) = f(qx) - f(x) \tag{3}$$

Consider the following expression:

$$\frac{f(x) - f(x_0)}{x - x_0} \tag{4}$$

As $x \to x_0$, the limit, if it exists, gives the definition of the derivative $\frac{df(x)}{dx}$ of a function $f(x)$ at $x=x_0$. In $q$-derivative, instead of independent variable $x$ of $f(x)$ be added by an infinitesimal value, it is multiplied by a $q$ as $x = qx_0$, $q \neq 1$ [37]. For function $f: R \to R$, and $q \neq 1$ the $q$-derivative can be expressed as follows:

$$D_q f(x) = \begin{cases} \frac{f(x) - f(qx)}{(1-q)(x)} & x \neq 0, q \neq 1 \\ \frac{\partial f(x)}{\partial x} & x = 0, q \neq 1 \end{cases} \quad (5)$$

For two functions $f(x)$ and $g(x)$ the $q$-differential and the $q$-derivative satisfy characteristics [6], [45], which are enlisted below.

1. The lack of symmetry in the differential of the product of two functions:

$$d_q(f(x)g(x)) = f(qx)d_q g(x) + g(x)d_q f(x) \quad (6)$$

2. For any constants $a$ and $b$ and two functions $f(x)$ and $g(x)$, the $q$-derivative of the sum of products can be determined as in (7):

$$D_q(af(x) + bg(x)) = aD_q f(x) + bD_q f(x) \quad (7)$$

3. The $q$-derivative of the product of $f(x)$ and $g(x)$ is given with the following formulae:

$$D_q(f(x)g(x)) = f(qx)D_q g(x) + g(x)D_q f(x) \quad (8)$$

$$D_q(f(x)g(x)) = f(x)D_q g(x) + g(qx)D_q f(x) \quad (9)$$

4. The $q$-derivative of the quotient of $f(x)$ and $g(x)$ is given as follows:

$$D_q \left(\frac{f(x)}{g(x)}\right) = \frac{g(x)D_q f(x) - f(x)D_q g(x)}{g(x)g(qx)} = \frac{g(qx)D_q f(x) - f(qx)D_q g(x)}{g(x)g(qx)} \quad (10)$$

In the limit $q \to 1$ the $q$-derivative tends to classical derivative [39], [44].

## 2.3 $q$-partial derivatives

For a real-valued continuous function of n variables $f(\mathbf{x})$, $f: R^n \to R$, the $q$-first-order partial derivative vector of the $n$ first-order partial derivatives with respect to a variable $x_i$ can be defined as follows [46], [47]:

$$D_{q_i,x_i} f(\mathbf{x}) = \frac{f(x_1, \cdots, x_{i-1}, x_i, x_{i+1}, \cdots, x_n) - f(x_1, \cdots, x_{i-1}, q_i x_i, x_{i+1}, \cdots, x_n)}{(1-q_i)x_i}, x_i \neq 0, q_i \neq 1 \quad (11)$$

Furthermore suppose the partial derivatives of $f_j: R^n \to R^m$, for $j=1,...,m$ exist. For $x \in R^n$ consider the operator $\varepsilon_{q,i}$ on $f_j$

$$(\varepsilon_{q,i} f_j)(x) = f_j(x_1, \cdots x_{i-1}, qx_i, x_{i+1}, \cdots x_n) \quad (12)$$

The $q$-partial derivative of $f_j$, for $j=1,...,m$ at $x$ with respect to $x_i$ for $i=1,...,n$ is then

$$D_{q,x_i} f_j(x) = \begin{cases} \frac{f_j(x) - \varepsilon_{q,i} f_j(x)}{(1-q)(x_i)} & x_i \neq 0 \\ \frac{\partial f_j(x)}{\partial x_i} & x_i = 0 \end{cases} \quad (13)$$

## 2.4 $q$-gradient

The concept of the $q$-gradient is derived from the definition of Jackson's derivative [48], as in (5), which is also known as the $q$-derivative [49]. Moreover, the concept of $q$-gradient vector is derived from the concept of gradient vector, also based on Jackson's derivative. For the parameter vector $q$ and function of $n$ variables $f(\mathbf{x})$, the $q$-gradient is given as follows:

$$\nabla_q f(\mathbf{x}) = [D_{q_1,x_1} f(\mathbf{x}) \cdots D_{q_i,x_i} f(\mathbf{x}) \cdots D_{q_n,x_n} f(\mathbf{x})] \quad (14)$$

In the limit $q_i \to 1$, for all $i$ ($i=1,...,n$) the $q$-gradient vector tends to classical gradient vector of partial derivatives of $f$ with respect to each dimension [4].

For any $g_i(\mathbf{x}) = \frac{\partial f_j(\mathbf{x})}{\partial x_i}$ ($i = 1, ..., n$ and $j = 1, ... m$)

$$g_j(\mathbf{x}) = \nabla f_j(\mathbf{x}) = \left[ g_1^j(\mathbf{x})\ g_2^j(\mathbf{x}) \cdots g_n^j(\mathbf{x}) \right]^T \tag{15}$$

The *q*-Jacobian matrix $\mathbf{J}_q(\mathbf{x})$ of the function $f_j$ for *j*=1,…,*m* represents matrix of *q*-partial derivatives of $\nabla f_j(\mathbf{x})$ [50]:

$$\mathbf{J}_q(\mathbf{x}) = D_q \nabla f_j(\mathbf{x}) = \begin{bmatrix} D_{q,x_1} g_1^j(x) & D_{q,x_2} g_1^j(x) & \cdots & D_{q,x_n} g_1^j(x) \\ D_{q,x_1} g_2^j(x) & D_{q,x_2} g_2^j(x) & \cdots & D_{q,x_n} g_2^j(x) \\ \vdots & \vdots & & \vdots \\ D_{q,x_1} g_n^j(x) & D_{q,x_2} g_n^j(x) & \cdots & D_{q,x_n} g_n^j(x) \end{bmatrix}_{nxn} \tag{16}$$

## 3 Gauss-Newton method

The GN method is a simple iterative procedure that resolves problems of convergence of the Newton's method in the non-linear least squares (NLS) optimization. Newton's method is a recursive algorithm for approximating the roots of differentiable function *f(x)* of single variable *x* that uses the first two terms of the Taylor's series of *f(x)* in the vicinity of a suspected root. The algorithm fails in certain cases leading to oscillations, divergence, or when the derivative is zero [51]. The GN method is a modification of Newton's method. In solving non-linear problems, the GN method is used to minimize the sum of quadratic function values, which does not require the calculation of the second-order derivatives [31]. Given **f**: $R^n \rightarrow R^m$, with *m≥n*, the GN method seeks to minimize the objective function

$$\text{minimize } f(\mathbf{x}) = \frac{1}{2} \sum_{i=1}^{m} r_i(\mathbf{x})^2 = \frac{1}{2} \mathbf{r}(\mathbf{x})^T \mathbf{r}(\mathbf{x}) = \frac{1}{2} \|\mathbf{r}(\mathbf{x})\|^2 \tag{17}$$

where $\mathbf{x} = [x_1\ x_2 \cdots x_n]^T$, $\mathbf{r} = [r_1(\mathbf{x}) \cdots r_m(\mathbf{x})]^T$, $\|\cdot\|$ is Euclidean norm and $\frac{1}{2}$ is used for convenience. The objective function **f(x)** is defined by *m* auxiliary residual functions $\{r_i(\mathbf{x})\}$. The GN approach to this optimization is to approximate **r** by a first order Taylor expansion in order to determine the search direction from the current iterate [52].

### 3.1 Geometric interpretation

Consider residuals $r_i(\mathbf{x})$ between model function and measured values. The residual function **r(x)** can be interpreted as a point in observation space $R^m$. Thus, the NLS problem can be interpreted as trying to find the point **x**$^*$ in parameter space $R^n$ that corresponds to the point **x**$^*$ in observation space $R^m$ closest to the origin. Consider the gradient of the function *f* given with [53]

$$\nabla \mathbf{f}(\mathbf{x}) = \nabla \mathbf{r}(\mathbf{x}) \mathbf{r}(\mathbf{x}) = \mathbf{J}(\mathbf{x})^T \mathbf{r}(\mathbf{x}) \tag{18}$$

$$\nabla f(\mathbf{x}) = \begin{bmatrix} \frac{\partial f(\mathbf{x})}{\partial x_1} \\ \vdots \\ \frac{\partial f(\mathbf{x})}{\partial x_n} \end{bmatrix} = \begin{bmatrix} r_1(\mathbf{x}) \frac{\partial r_1(\mathbf{x})}{\partial x_1} + r_2(\mathbf{x}) \frac{\partial r_2(\mathbf{x})}{\partial x_1} + \cdots + r_m(\mathbf{x}) \frac{\partial r_m(\mathbf{x})}{\partial x_1} \\ \vdots \\ r_1(\mathbf{x}) \frac{\partial r_1(\mathbf{x})}{\partial x_n} + r_2(\mathbf{x}) \frac{\partial r_2(\mathbf{x})}{\partial x_n} + \cdots + r_m(\mathbf{x}) \frac{\partial r_m(\mathbf{x})}{\partial x_n} \end{bmatrix} \tag{19}$$

The Jacobian matrix **J(x)** of the first partial derivatives and the Hessian matrix of second partial derivatives **H**, are given in (20) and (21), respectively.

$$\mathbf{J}(\mathbf{x}) = \begin{bmatrix} \frac{\partial r_1(\mathbf{x})}{\partial x_1} & \frac{\partial r_2(\mathbf{x})}{\partial x_1} & \cdots & \frac{\partial r_m(\mathbf{x})}{\partial x_1} \\ & & \vdots & \\ \frac{\partial r_1(\mathbf{x})}{\partial x_n} & \frac{\partial r_2(\mathbf{x})}{\partial x_n} & \cdots & \frac{\partial r_m(\mathbf{x})}{\partial x_n} \end{bmatrix} \tag{20}$$

$$\mathbf{H} = \nabla^2 f(\mathbf{x}) = \nabla \mathbf{r}(\mathbf{x}) \mathbf{r}(\mathbf{x}) + \sum_{i=1}^{m} r_i(\mathbf{x}) \nabla^2 r_i(\mathbf{x}) = \mathbf{J}(\mathbf{x})^T \mathbf{J}(\mathbf{x}) + \mathbf{Q}(\mathbf{x}) \tag{21}$$

The Hessian matrix consists of two terms: $\mathbf{J}(\mathbf{x})^T \mathbf{J}(\mathbf{x})$ with only the first derivatives and $\mathbf{Q}(\mathbf{x})$ with second order derivatives.

The coefficients of the Hessian matrix at the position $k$, $j$ are determined as

$$H_{k,j} = \frac{\partial^2 \mathbf{f}(\mathbf{x})}{\partial x_k \partial x_j} = \sum_{i=1}^{m} r_i(\mathbf{x}) \frac{\partial^2 r_i(\mathbf{x})}{\partial x_k \partial x_j} \tag{22}$$

where $i = 1, \dots, m$ and $j, k = 1, \dots, n$. Size of the Hessian matrix is one of the main challenges in Newton's and Newton-like algorithms. The Newton's method suffers from lack of robustness and might not converge since it is not guaranteed that it fulfills the descent condition for NLS problems. The GN method improves robustness by simplification of the Newton's method with $\mathbf{J}(\mathbf{x})^T \mathbf{J}(\mathbf{x})$ approximation of the Hessian matrix near a given point [54]. Method determines the search direction $\mathbf{h}_{GN}$ as the solution of the Newton's equation [55]:

$$\nabla^2 \mathbf{f}(\mathbf{x}) \mathbf{h}_N = -\nabla \mathbf{f}(\mathbf{x}) \tag{23}$$

where $\mathbf{h}_N$ represents the search direction of the Newton's method. Accordingly, the search direction can be determined as [56]:

$$\mathbf{J}(\mathbf{x})^T \mathbf{J}(\mathbf{x}) \mathbf{h}_{GN} = -\mathbf{J}(\mathbf{x})^T \mathbf{f}(\mathbf{x}) \tag{24}$$

The GN algorithm is particularly suited to the small residual case. Therefore, the search direction $\mathbf{h}_{GN}$ can be found as the solution of the following equation:

$$\mathbf{h}_{GN} = -\left(\mathbf{J}(\mathbf{x})^T \mathbf{J}(\mathbf{x})\right)^{-1} \mathbf{J}(\mathbf{x})^T f(\mathbf{x}) = -\mathbf{J}(\mathbf{x})^\dagger \mathbf{f}(\mathbf{x}) \tag{25}$$

where $\mathbf{J}(\mathbf{x})^\dagger = \left(\mathbf{J}(\mathbf{x})^T \mathbf{J}(\mathbf{x})\right)^{-1} \mathbf{J}(\mathbf{x})^T$ denotes Moore-Penrose pseudo-inverse of $\mathbf{J}(\mathbf{x})$ [57], [58], [59]. If $\mathbf{J}(\mathbf{x})$ has a full rank, the Hessian approximation $\mathbf{J}(\mathbf{x})^T \mathbf{J}(\mathbf{x})$, is positive definite and the Gauss-Newton search direction is a descent direction [60], [61] where

$$f(\mathbf{x}^*) \leq f(\mathbf{x}), \ \|\mathbf{x} - \mathbf{x}^*\| < \delta \tag{26}$$

for a small positive number $\delta$. Otherwise $\mathbf{J}(\mathbf{x})^T \mathbf{J}(\mathbf{x})$ is non-invertible and the equation (24) does not have unique solution. In this case, problem is said to be under-determined or over-parameterized. The advantage of GN over the Newton's method is that the second order derivatives $\nabla^2 r_i(\mathbf{x})$ need not to be calculated. If $\mathbf{r}(\mathbf{x}^*) = 0$ then the approximation $\mathbf{Q}(\mathbf{x}) \approx 0$ is good enough and the GN method behave like the Newton's method close to the solution, i.e. the method has quadratic convergence if $\mathbf{J}(\mathbf{x}^*)$ has full rank. However, if any residual component $r_i(\mathbf{x}^*)$ and corresponding curvature $\nabla^2 r_i(\mathbf{x})$ is large, the GN method will converge slower that the Newton's method or may not even be locally convergent. With local strategy such as line search, it would converge independently of how close to the solution algorithm starts. The initial guess of GN method is $\mathbf{x}=\mathbf{x}^{(0)}$ [62]. Then $\mathbf{x}^{(k+1)}$ can be calculated iteratively, as follows:

$$\mathbf{x}^{(k+1)} = \mathbf{x}^{(k)} + \alpha^{(k)} \mathbf{h}_{GN}^{(k)} \tag{27}$$

where $\alpha^{(k)} \in (0, 1]$ improves the convergence and can be found by the line search [63]. Usually, $\alpha^{(k)}$ is set to one. The method for line search has guaranteed convergence provided that $\{\mathbf{x} | f(\mathbf{x}) \leq f(\mathbf{x}_0)\}$ is bounded and that Jacobian $\mathbf{J}(\mathbf{x})$ has full rank in all steps.

### 3.2 Analytical approach

The analytical approach considers the Taylor expansion of the function $f: \mathbb{R}^n \to \mathbb{R}^m$, with $m \geq n$. Minimize $\|f(\mathbf{x})\|$ so that [64]

$$\mathbf{x}^* = \arg\min_{\mathbf{x}}\{F(\mathbf{x})\} \tag{28}$$

where $F(\mathbf{x})$ is the SSE given with (29)

$$F(\mathbf{x}) = \frac{1}{2}\sum_{i=1}^{m}\left(f_i(\mathbf{x})\right)^2 = \frac{1}{2}\|f(\mathbf{x})\|^2 = \frac{1}{2}\mathbf{f}(\mathbf{x})^T \mathbf{f}(\mathbf{x}) \tag{29}$$

Provided that $\mathbf{f}$ has continuous second partial derivatives then [55]:

$$\mathbf{f}(\mathbf{x} + \mathbf{h}) = \mathbf{f}(\mathbf{x}) + \mathbf{J}(\mathbf{x})\mathbf{h} + O(\|\mathbf{h}\|^2) \tag{30}$$

where **h** is the step and

$$\mathbf{J(x)} \in \mathbf{R}^{m \times n} \tag{31}$$

represents Jacobian matrix containing the first partial derivatives of the components of function **f**. At the position (*i,j*) (*i*=1,…,*n*; *j*=1,…,*m*), elements of the Jacobian matrix are given with the formula

$$[\mathbf{J(x)}]_{ij} = \frac{\partial f_i(\mathbf{x})}{\partial x_j} \tag{32}$$

Therefore, the SSE can be expressed as

$$\frac{\partial F(\mathbf{x})}{\partial x_j} = \sum_{i=1}^{m} f_i(\mathbf{x}) \frac{\partial f_i(\mathbf{x})}{\partial x_j} = f_1(\mathbf{x}) \frac{\partial f_1}{x_j} + \cdots + f_m(\mathbf{x}) \frac{\partial f_m}{x_j} \tag{33}$$

The GN method is based on a linear approximation to the components of **f** in the neighborhood of **x**. For small ‖**h**‖ the Taylor expansion becomes can be approximated as follows [55]

$$\mathbf{f(x+h)} \approx l(\mathbf{h}) \equiv \mathbf{f(x)} + \mathbf{J(x)h} \tag{34}$$

which gives the following formulae

$$F(\mathbf{x+h}) \approx L(\mathbf{h}) \equiv \tfrac{1}{2} l(\mathbf{h})^\mathrm{T} l(\mathbf{h}) = \tfrac{1}{2} \mathbf{f}^\mathrm{T} \mathbf{f} + \mathbf{h}^\mathrm{T} \mathbf{J}^\mathrm{T} \mathbf{f} + \tfrac{1}{2} \mathbf{h}^\mathrm{T} \mathbf{J}^\mathrm{T} \mathbf{J} \mathbf{h} \tag{35}$$

$$F(\mathbf{x+h}) = F(\mathbf{x}) + \mathbf{h}^\mathrm{T} \mathbf{J}^\mathrm{T} \mathbf{f} + \tfrac{1}{2} \mathbf{h}^\mathrm{T} \mathbf{J}^\mathrm{T} \mathbf{J} \mathbf{h} \tag{36}$$

with **f**=**f**(**x**), **J**=**J**(**x**) and **h**=**h**$_{GN}$. Thus, the matrix of the gradients is given with [65]

$$\mathbf{F'(x)} = \mathbf{J(x)}^\mathrm{T} \mathbf{f(x)} \tag{37}$$

showing that [55]

$$\mathbf{F''(x)} = \mathbf{J(x)}^\mathrm{T} \mathbf{f(x)} + \sum_{i=1}^{m} f_i(\mathbf{x}) \mathbf{f}_i''(\mathbf{x}) \tag{38}$$

Since the GN method is based on linear approximation of components of **f** in the neighborhood of **x** and ‖**h**‖ suppose to be a small number, the Taylor expansion can be approximated as follows [55]

$$\mathbf{f(x+h)} \cong l(\mathbf{h}) = \mathbf{f(x)} + \mathbf{J(x)h} \tag{39}$$

and

$$F(\mathbf{x+h}) \cong L(\mathbf{h}) = \tfrac{1}{2} l(\mathbf{h})^\mathrm{T} l(\mathbf{h}) = F(\mathbf{x}) + \mathbf{h}^\mathrm{T} \mathbf{J}^\mathrm{T} \mathbf{f} + \tfrac{1}{2} \mathbf{h}^\mathrm{T} \mathbf{J}^\mathrm{T} \mathbf{J} \mathbf{h} \tag{40}$$

with **f** = **f**(**x**) and **J** = **J**(**x**). The GN step **h**$_{GN}$ minimizes *L*(**h**) so that

$$\mathbf{h}_{GN} = \mathrm{argmin}_\mathbf{h}\{L(\mathbf{h})\} \tag{41}$$

$$\mathbf{L'(h)} = \mathbf{J}^\mathrm{T}\mathbf{f} + \mathbf{J}^\mathrm{T}\mathbf{J}\mathbf{h}, \quad \mathbf{L''(h)} = \mathbf{J}^\mathrm{T}\mathbf{J} \tag{42}$$

where **L'**(**0**)=**F'**(**x**) and **L''**(**h**) is symmetric and independent of **h**. If **J** has full rank (the columns are linearly independent) than **L''**(**h**) is positive definite, which implies that **L**(**h**) has unique minimizer that can be found by solving the equation

$$(\mathbf{J}^\mathrm{T}\mathbf{J})\mathbf{h}_{GN} = -\mathbf{J}^\mathrm{T}\mathbf{f} \tag{43}$$

This is a descent direction of **F** [55] since

$$\mathbf{h}_{GN}^\mathrm{T} \mathbf{F'(x)} = \mathbf{h}_{GN}^\mathrm{T} (\mathbf{J}^\mathrm{T}\mathbf{J}) \mathbf{h}_{GN} < 0 \tag{44}$$

The initial guess and the search direction is the same as in (27). To compare the search directions $\mathbf{h}_N$ and $\mathbf{h}_{GN}$ consider the following formulae

$$\mathbf{F}''(\mathbf{x})\mathbf{h}_N = -\mathbf{F}'(\mathbf{x}), \ \mathbf{L}''(\mathbf{h})\mathbf{h}_{GN} = -\mathbf{L}'(\mathbf{0}) \tag{45}$$

Since the right-hand sides of the formulae are identical, then

$$\mathbf{F}''(\mathbf{x}) = \mathbf{L}''(\mathbf{h}) + \sum_{i=1}^{m} f_i(\mathbf{x})\mathbf{f}_i''(\mathbf{x}) \tag{46}$$

If $\mathbf{f}(\mathbf{x}^*)=\mathbf{0}$ then $\mathbf{L}''(\mathbf{h}) \cong \mathbf{F}''(\mathbf{x})$ for $\mathbf{x}$ close to $\mathbf{x}^*$ and GN algorithm gets quadratic convergence. Note that $F(\mathbf{x}^*)$ controls the convergence speed. It is well known that GN algorithm may fail to be well defined [66], [67]. To ensure that the GN method is well defined and converges to the stationary point, the initial iterate has to be close enough to the solution. In [67] authors have given an estimate of convergence radius and a criterion for convergence of the GN method [68].

## 4  q-Gauss-Newton method

The q-GN method as well as the classical GN method can only be used to minimize sum of squared function values and has the advantages that the q-Hessian matrix is not required. The q-GN method can be derived by linearly approximating the vector of functions $f_i$ using q-analog to the Taylor's formula. In [44], authors derived q-Taylor's formula for any polynomial $f(x)$ of degree $n$ and any number $c$. The truncated q-Taylor expansion can be expressed as follows

$$f(x) = \sum_{j=0}^{n}(D_q^j f)(c)\frac{(x-c)_q^j}{j!} \tag{47}$$

$$f(x) = 1 + (D_q f)(c)(x-c)_q + (D_q^2 f)(c)\frac{(x-c)_q^2}{2!} + \cdots + (D_q^n f)(c)\frac{(x-c)_q^n}{n!} \tag{48}$$

where $D_q^j f$ represents $j$-th order q-derivative of the objective function $f$. Suppose that $\mathbf{f}(\mathbf{x}^*)=\mathbf{0}$. The q-Taylor series of the $\mathbf{f}(\mathbf{x})$ around $\mathbf{x}^*$ can be approximated with [69]

$$\mathbf{f}(\mathbf{x}^*) \approx \mathbf{f}(\mathbf{x}) + D_{q,\mathbf{x}}\mathbf{f}(\mathbf{x})(\mathbf{x}^* - \mathbf{x}) \tag{49}$$

where $D_{q,\mathbf{x}}\mathbf{f}(\mathbf{x}) = \left[D_{q,x_j}f_i(\mathbf{x})\right]_{mxn}$ represents the q-Jacobian matrix of q-partial derivatives of $\mathbf{f}$. Considering (12) and (39) the q-Taylor expansion of $\mathbf{f}$ can be expressed as

$$\mathbf{f}(\mathbf{x}+\mathbf{h}) \approx \mathbf{f}(\mathbf{x}) + J_q(\mathbf{x})\mathbf{h} \tag{50}$$

For small value of $\mathbf{h}$, models $l_q(\mathbf{h})$ and $L_q(\mathbf{h})$ can be approximated as in (39) and (40), respectively, so that

$$\mathbf{f}(\mathbf{x}+\mathbf{h}) \cong l_q(\mathbf{h}) = \mathbf{f}(\mathbf{x}) + J_q(\mathbf{x})\mathbf{h} \tag{51}$$

$$F(\mathbf{x}+\mathbf{h}) \cong L_q(\mathbf{h}) = F(\mathbf{x}) + \mathbf{h}^T J_q^T \mathbf{f} + \frac{1}{2}\mathbf{h}^T J_q^T J_q \mathbf{h} \tag{52}$$

where

$$J_q = D_{q,x_j}F(\mathbf{x}) = \sum_{i=1}^{m}f_i(\mathbf{x})D_{q,x_j}f_i(\mathbf{x}) = \sum_{i=1}^{m}f_i(\mathbf{x})\frac{f_i(x_1,\cdots,x_{j-1},x_j,x_{j+1},\cdots x_n)-f_i(x_1,\cdots,x_{j-1},qx_j,x_{j+1},\cdots x_n)}{(1-q)x_j} \tag{53}$$

If $J_q$ has full rank than $L_q(\mathbf{h})$ has minimizer which can be found by solving the equation

$$(J_q^T J_q)\mathbf{h}_{q-GN} = -J_q^T \mathbf{f} \tag{54}$$

A descent direction for $F$ if the q-GN step $\mathbf{h}_{q-GN}$ satisfies the relation:

$$\mathbf{h}_{q-GN}{}^T \mathbf{F}'(\mathbf{x}) = \mathbf{h}_{q-GN}{}^T (J_q^T J_q)\mathbf{h}_{q-GN} < 0 \tag{55}$$

where $\mathbf{F}'(\mathbf{x}) = J_q(\mathbf{x})^T \mathbf{f}(\mathbf{x})$.

Therefore

$$\mathbf{h}_{q-\text{GN}} = -\left(\mathbf{J}_q(\mathbf{x})^{\text{T}} \mathbf{J}_q(\mathbf{x})\right)^{-1} \mathbf{J}_q(\mathbf{x})^{\text{T}} \mathbf{f}(\mathbf{x}) = -\mathbf{J}_q(\mathbf{x})^{\dagger} \mathbf{f}(\mathbf{x}) \qquad (56)$$

where $\mathbf{J}_q(\mathbf{x})^{\dagger} = \left(\mathbf{J}_q(\mathbf{x})^{\text{T}} \mathbf{J}_q(\mathbf{x})\right)^{-1} \mathbf{J}_q(\mathbf{x})^{\text{T}}$ denotes Moore-Penrose pseudo-inverse of $\mathbf{J}_q(\mathbf{x})$. Thus $\mathbf{x}^{(k+1)}$ can be evaluated as:

$$\mathbf{x}^{(k+1)} = \mathbf{x}^{(k)} + \alpha^{(k)} \mathbf{h}_{q-\text{GN}}^{(k)} \qquad (57)$$

at each iteration step of the *q*-GN algorithm. Considering processing power, calculating the *q*-Jacobian matrix is less demanding then calculation of the *q*-Hessian matrix. In this way the *q*-GN algorithm can be faster than Newton's algorithm. However, the *q*-Jacobian matrix tends to be inaccurate when the initial point is set far from the solution since it is approximation based on the assumption that the residuals are small, which may cause the *q*-GN algorithm to fail [70]. The *q*-GN algorithm can be expressed in pseudocode as:

1. analytically determine the $Jq$ of $f$
2. set:
   {counter $k = 0$;
   initial guess $x(0)$;
   stop iteration $stop\_iter$ to small positive value;
   maximum number of iterations $max\_no\_iter$;
   $q$;}
3. while(1)
   evaluate $f(k)$;
   evaluate $Jq(k)$;
   evaluate $\Box q(k)$;
   determine $x(k+1) = x(k) + \Box q(k)$;
   calculate $norm(k)$;
   if $(norm(k) \leq stop\_iter)$ or $(k = \max\_no\_iter)$
       break
   end
end

Since the *q*-GN algorithm requires calculation of the *q*-Jacobian matrix at each iteration step the general formula for the *q*-Jacobian should be analytically determined before the *q*-GN algorithm starts. It should be noted that choosing *q*-parameter to be close to one may cause the algorithm to diverge or oscillate. This happens in some cases when the objective function contains sinusoidal-based terms, exponential values or absolutes values of unknown variables.

## 5 Experimental results

The convergence of *q*-GN method is assured only when the initial guess $x_0$ is close enough to the solution. The calculation starts from an initial guess and at each iteration step the approximate solution $x_n$ is calculated. The iteration stops when the convergence is achieved within a given tolerance. In general, the initial guess for the solution may be outside the region of convergence. An important strategy is that when $x_n$ approaches to the solution the convergence rate increases. The exit/stopping criteria may be the minimum SSE or predefined maximum number of iterations, which is the safeguard against infinite loops.

In this research we have evaluated the performances of the *q*-GN method and observed the change of the algorithm response with the increasing of the *q*-parameter for three examples. The convergences of the proposed algorithms are examined for $q \in \{0.9, 0.95, 0.99, 0.9995, 1\}$.

**Example 1**: In this example a global minimum of the function $f(x) = 2 - \left(e^{-x^2} + 2e^{-(x-3)^2}\right)$ is found (See Figure 1). The $q$-Jacobian of the function is given with $J_q = \frac{-e^{-x^2} + e^{-(qx)^2} - 2e^{-(x-3)^2} + 2e^{-(qx-3)^2}}{(1-q)x}, q \neq 1, x \neq 0$. The initial guess $x_0 = 2.1$ is chosen to be near the optimal solution. The global minimum of the function is $f(x_{global}) = 0$, for $x_{global} = 3$.

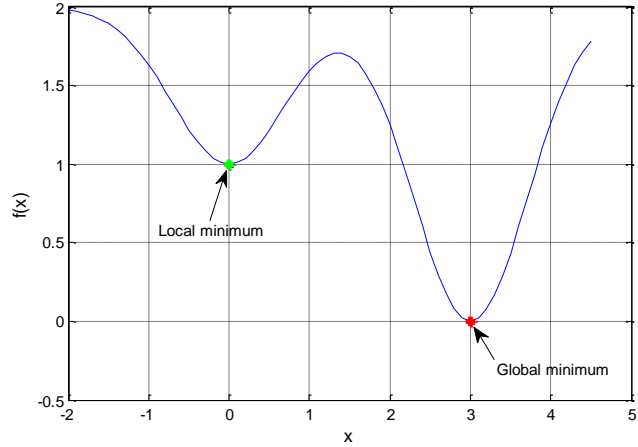

Figure 1: Function $f(x) = 2 - \left(e^{-x^2} + 2e^{-(x-3)^2}\right)$

Table 1 shows the results on *x* and *f*(*x*), for *q*=0.9, *q*=0.95, *q*=0.99 and *q*=0.9995. The number of iterations decreases as the *q* parameter decreases from 19 down to 8, which points to the high-speed calculations when *q*→1. The norm of the objective function also decreases with the increase of the *q*-parameter. Results show that the *q*-GN method decreases the number of calculations of the iterative process near the global minimum. Significant decreasing of the number of iterations as the parameter *q*→1 shows that the *q*-parameter can be considered an additional fitting factor.

Table 1: *x* and *f*(*x*) for *q*∈{0.9, 0.95, 0.99, 0.9995}

| q=0.9995 Iterations = 8 norm=63131·10⁻⁶ | | q=0.99 Iterations = 10 norm=4.6782·10⁻⁵ | | q=0.95 Iterations = 15 norm=2.3923·10⁻⁴ | | q=0.9 Iterations = 19 norm=5.0199·10⁻⁴ | |
|---|---|---|---|---|---|---|---|
| x | f(x) | x | f(x) | x | f(x) | x | f(x) |
| 2.8168 | 0.2946 | 2.8125 | 0.2946 | 2.7984 | 0.2946 | 2.7879 | 0.2946 |
| 2.9094 | 0.0657 | 2.9015 | 0.0687 | 2.8776 | 0.0793 | 2.8599 | 0.0875 |
| 2.9541 | 0.0162 | 2.9444 | 0.0191 | 2.9173 | 0.0295 | 2.8977 | 0.0386 |
| 2.9759 | 0.0041 | 2.9659 | 0.0060 | 2.9396 | 0.0134 | 2.9204 | 0.0206 |
| 2.9863 | 0.0010 | 2.9771 | 0.0022 | 2.9532 | 0.0071 | 2.9352 | 0.0124 |
| 2.9905 | 0.0002 | 2.9832 | 0.0009 | 2.9621 | 0.0042 | 2.9456 | 0.0082 |
| 2.9917 | 0.0000 | 2.9866 | 0.0004 | 2.9683 | 0.0027 | 2.9532 | 0.0057 |
| 2.9919 | 0.0000 | 2.9886 | 0.0002 | 2.9728 | 0.0019 | 2.9590 | 0.0042 |
| | | 2.9898 | 0.0001 | 2.9761 | 0.0013 | 2.9635 | 0.0032 |
| | | 2.9906 | 0.0001 | 2.9787 | 0.0010 | 2.9671 | 0.0025 |
| | | | | 2.9808 | 0.0008 | 2.9700 | 0.0020 |
| | | | | 2.9824 | 0.0006 | 2.9725 | 0.0016 |
| | | | | 2.9837 | 0.0005 | 2.9745 | 0.0014 |
| | | | | 2.9848 | 0.0004 | 2.9763 | 0.0012 |
| | | | | 2.9857 | 0.0003 | 2.9778 | 0.0010 |
| | | | | | | 2.9791 | 0.0008 |
| | | | | | | 2.9802 | 0.0007 |
| | | | | | | 2.9812 | 0.0006 |

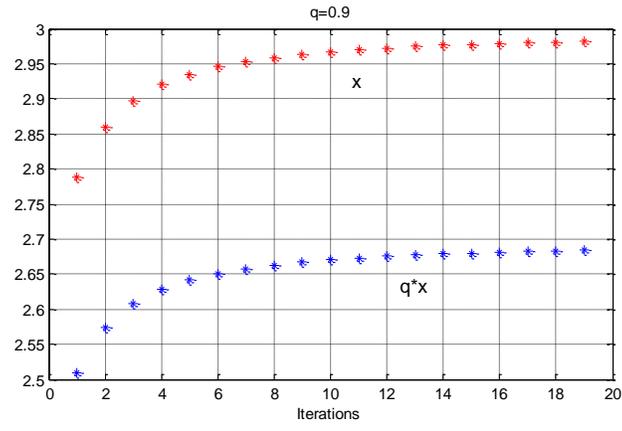

Figure 2: *x* and *q*x* (Iterations=19, *q*=0.9)

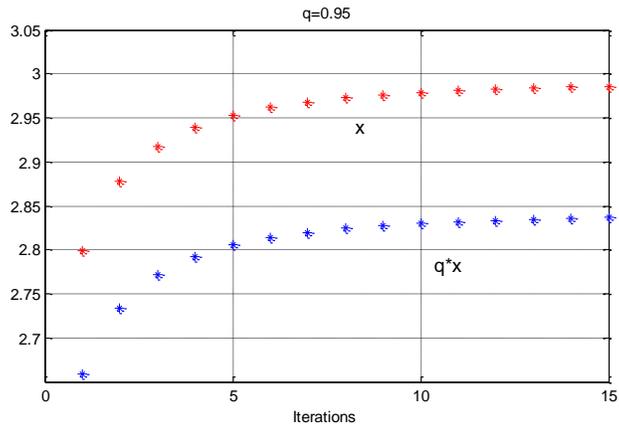

Figure 3: *x* and *q*x* (Iterations=15, *q*=0.95)

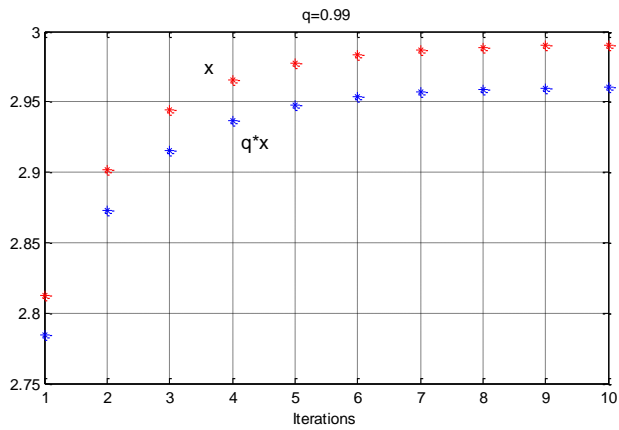

Figure 4: *x* and *q*x* (Iterations=10, *q*=0.99)

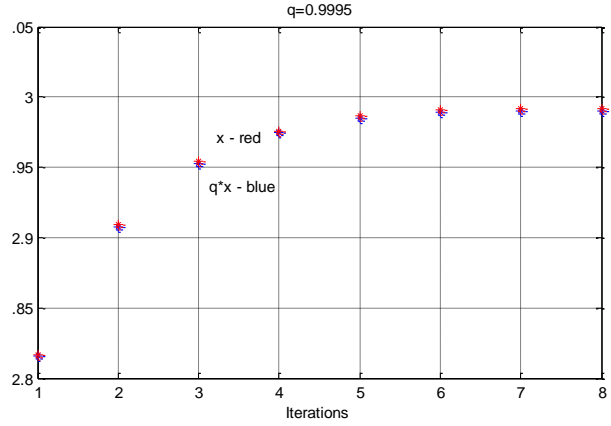

Figure 5: *x* and *q\*x* (Iterations=8, *q*=0.9995)

We compared the results with a Nelder-Mead direct search method, described in [42] and designed in MATLAB as a function for finding minimum of unconstrained multivariable using derivative-free method. Although estimated values of $x_{optim} = 2.9998$ and $f(x_{optim}) = -1.2348 \cdot 10^{-4}$ are more precise for the Nelder-Mead search, the number of iterations is 16. Obviously, the *q*-GN algorithm tends to diminishes number of calculations to reach the optimization goal, as the parameter *q*→1.

***Example 2***: This example considers the problem taken from Powell (1970) [71], and described in [55]. For a real-valued function of two variables ***f*(x)**, ***f***: R²→R², $f(x_1, x_2) = \begin{bmatrix} f_1 \\ f_2 \end{bmatrix} = \begin{bmatrix} x_1 \\ \frac{10x_1}{x_1+0.1} + 2x_2^2 \end{bmatrix}, x_1 \neq -0.1$, find optimal solution, with **x*** =**0** only solution and Jacobian $J = \begin{bmatrix} 1 & 0 \\ (x_1+0.1)^{-2} & 4x_2 \end{bmatrix}$, $x_1 \neq -0.1$ singular at the solution. Authors took $x_0 = [3\ 1]^T$ and applied algorithm with exact line search. The iterates converged and the iteration process stopped for $x_c = [1.8016\ 0]^T$, which is not a solution. On the other hand, better solution is determined with the GN algorithm. The algorithm stopped for $x_c = [3.82 \cdot 10^{-8}\ -1.38 \cdot 10^{-3}]^T$. The *q*-GN algorithm, with *q*-Jacobian $J_q = \begin{bmatrix} 1 & 0 \\ \frac{1}{(x_1+0.1)(qx_1+0.1)} & 2(1+q)x_2 \end{bmatrix}$, $x_1 \neq -0.1$, $q \neq 1$, determined the unique optimal solution which is $f(x_{optim}) = \begin{bmatrix} 0 \\ 0 \end{bmatrix}$, for $x_{optim} = [0\ 0]^T$. The iterative procedure starts with $x_0 = [-1\ 1]^T$, close to the solution. The results related to the number of iterations are presented in Table 2.

Table 2: $x_1$, $x_2$, $f_1$ and $f_2$, for *q*∈{0.9, 0.95, 0.99, 0.9995}

| q=0.9995 Iterations = 19 norm=9.8311·10⁻⁷ | | | | q=0.99 Iterations = 19 norm=9.1778·10⁻⁶ | | | | q=0.95 Iterations = 19 norm=6.7851·10⁻⁶ | | | | q=0.9 Iterations = 18 norm=9.5207·10⁻⁶ | | | |
|---|---|---|---|---|---|---|---|---|---|---|---|---|---|---|---|
| x₁ | x₂ | f₁ | f₂ | x₁ | x₂ | f₁ | f₂ | x₁ | x₂ | f₁ | f₂ | x₁ | x₂ | f₁ | f₂ |
| 0.0000 | -2.8575 | -1.0000 | 13.1111 | -0.0025 | -2.6970 | -1.0000 | 13.1111 | 0.0000 | -2.6970 | -1.0000 | 13.1111 | 0.0000 | -2.8158 | -1.0000 | 13.1111 |
| 0.0000 | -1.2934 | 0.0000 | 13.3902 | 0.0000 | -1.2974 | 0.0000 | 13.6026 | 0.0000 | -1.3139 | 0.0000 | 14.5476 | 0.0000 | -1.3338 | 0.0000 | 15.8573 |
| 0.0000 | -0.6465 | 0.0000 | 3.3459 | 0.0000 | -0.6454 | 0.0000 | 3.3666 | 0.0000 | -0.6401 | 0.0000 | 3.4528 | 0.0000 | -0.6318 | 0.0000 | 3.5580 |
| 0.0000 | -0.3232 | 0.0000 | 0.8361 | 0.0000 | -0.3211 | 0.0000 | 0.8332 | 0.0000 | -0.3199 | 0.0000 | 0.8195 | 0.0000 | -0.2993 | 0.0000 | 0.7983 |
| 0.0000 | -0.1616 | 0.0000 | 0.2089 | 0.0000 | -0.1597 | 0.0000 | 0.2062 | 0.0000 | -0.1519 | 0.0000 | 0.1945 | 0.0000 | -0.1418 | 0.0000 | 0.1791 |
| 0.0000 | -0.0808 | 0.0000 | 0.0522 | 0.0000 | -0.0795 | 0.0000 | 0.0510 | 0.0000 | -0.0740 | 0.0000 | 0.0462 | 0.0000 | -0.0671 | 0.0000 | 0.0402 |
| 0.0000 | -0.0404 | 0.0000 | 0.0130 | 0.0000 | -0.0395 | 0.0000 | 0.0126 | 0.0000 | -0.0361 | 0.0000 | 0.0110 | 0.0000 | -0.0318 | 0.0000 | 0.0090 |
| 0.0000 | -0.0202 | 0.0000 | 0.0033 | 0.0000 | -0.0197 | 0.0000 | 0.0031 | 0.0000 | -0.0176 | 0.0000 | 0.0026 | 0.0000 | -0.0151 | 0.0000 | 0.0020 |
| 0.0000 | -0.0101 | 0.0000 | 0.0008 | 0.0000 | -0.0098 | 0.0000 | 0.0008 | 0.0000 | -0.0086 | 0.0000 | 0.0006 | 0.0000 | -0.0071 | 0.0000 | 0.0005 |
| 0.0000 | -0.0050 | 0.0000 | 0.0002 | 0.0000 | -0.0049 | 0.0000 | 0.0002 | 0.0000 | -0.0042 | 0.0000 | 0.0001 | 0.0000 | -0.0034 | 0.0000 | 0.0001 |
| 0.0000 | -0.0025 | 0.0000 | 0.0001 | 0.0000 | -0.0024 | 0.0000 | 0.0000 | 0.0000 | -0.0020 | 0.0000 | 0.0000 | 0.0000 | -0.0016 | 0.0000 | 0.0000 |
| 0.0000 | -0.0013 | 0.0000 | 0.0000 | 0.0000 | -0.0012 | 0.0000 | 0.0000 | 0.0000 | -0.0010 | 0.0000 | 0.0000 | 0.0000 | -0.0008 | 0.0000 | 0.0000 |
| 0.0000 | -0.0006 | 0.0000 | 0.0000 | 0.0000 | -0.0006 | 0.0000 | 0.0000 | 0.0000 | -0.0005 | 0.0000 | 0.0000 | 0.0000 | -0.0004 | 0.0000 | 0.0000 |
| 0.0000 | -0.0003 | 0.0000 | 0.0000 | 0.0000 | -0.0003 | 0.0000 | 0.0000 | 0.0000 | -0.0002 | 0.0000 | 0.0000 | 0.0000 | -0.0002 | 0.0000 | 0.0000 |
| 0.0000 | -0.0002 | 0.0000 | 0.0000 | 0.0000 | -0.0001 | 0.0000 | 0.0000 | 0.0000 | -0.0001 | 0.0000 | 0.0000 | 0.0000 | -0.0001 | 0.0000 | 0.0000 |
| 0.0000 | -0.0001 | 0.0000 | 0.0000 | 0.0000 | -0.0001 | 0.0000 | 0.0000 | 0.0000 | 0.0000 | 0.0000 | 0.0000 | 0.0000 | 0.0000 | 0.0000 | 0.0000 |
| 0.0000 | -0.0001 | 0.0000 | 0.0000 | 0.0000 | 0.0000 | 0.0000 | 0.0000 | 0.0000 | 0.0000 | 0.0000 | 0.0000 | 0.0000 | 0.0000 | 0.0000 | 0.0000 |
| 0.0000 | -0.0001 | 0.0000 | 0.0000 | 0.0000 | 0.0000 | 0.0000 | 0.0000 | 0.0000 | 0.0000 | 0.0000 | 0.0000 | 0.0000 | 0.0000 | 0.0000 | 0.0000 |
| 0.0000 | 0.0000 | 0.0000 | 0.0000 | 0.0000 | 0.0000 | 0.0000 | 0.0000 | 0.0000 | 0.0000 | 0.0000 | 0.0000 | | | | |

Although we have assumed that the parameter $q$ improved the performances of the calculations Table 2 shows that the number of iterations negligible varies from 18 to 19, regardless of the value of $q$. Obviously, the $q$-Jacobian does not decrease the number of iterations. However, the results show that $q$-GN algorithm gives more accurate results than those presented in [36]. It should be noted that all optimization algorithms are designed and implemented as MATLAB functions. Figures 6-9 show the results on fitting $x_1$ and $x_2$ during iterative processes, in relation to the parameter $q$.

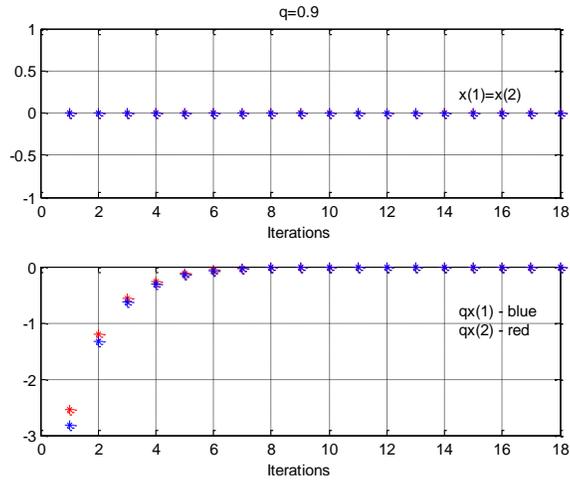

Figure 6: $x_1$ and $x_2$ (upwards) and $q^* x_1$ and $q^* x_2$ (downwards) (Iterations=18, $q$=0.9)

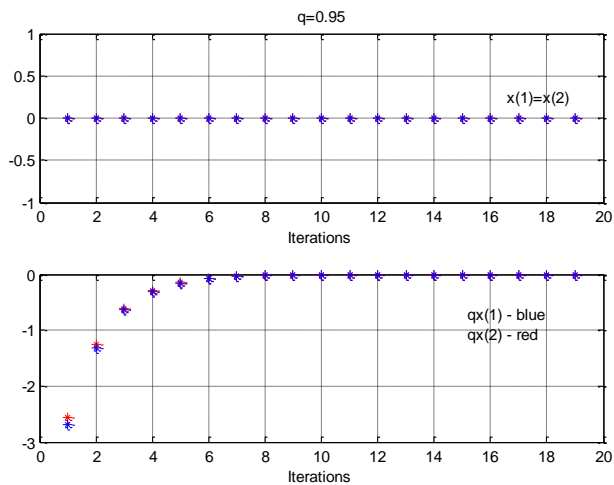

Figure 7: $x_1$ and $x_2$ (upwards) and $q^* x_1$ and $q^* x_2$ (downwards) (Iterations=19, $q$=0.95)

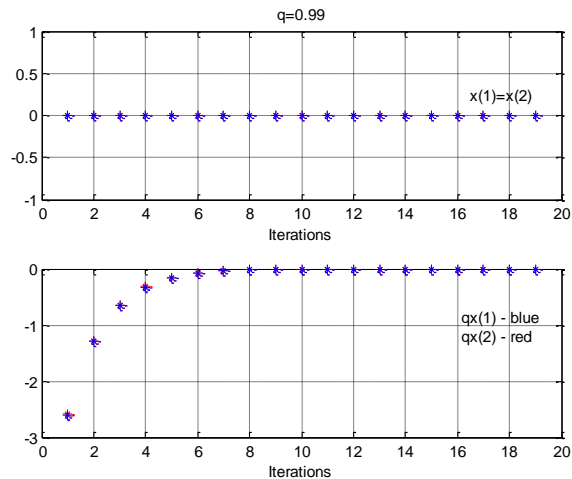

Figure 8: $x_1$ and $x_2$ (upwards) and $q^* x_1$ and $q^* x_2$ (downwards) (Iterations=19, $q$=0.99)

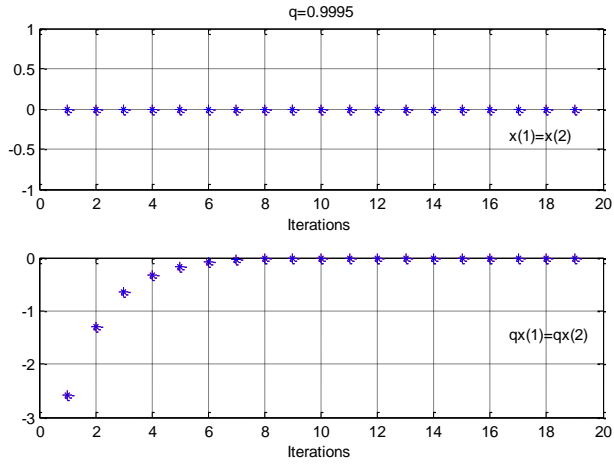

Figure 9: $x_1$ and $x_2$ (upwards) and $q*x_1$ and $q*x_2$ (downwards) (Iterations=19, $q$=0.9995)

**Example 3**: For a real-valued function of two variables $f(x)$, $f: R^2 \rightarrow R^3$, given with $f(x_1, x_2) = \begin{bmatrix} f_1 \\ f_2 \\ f_3 \end{bmatrix} = \begin{bmatrix} x_1 - 0.4 \\ x_2 - 8 \\ x_1^2 + x_2^2 - 1 \end{bmatrix}$, the $q$-Jacobian can be determined with the formula $J_q = \begin{bmatrix} 1 & 0 \\ 0 & 1 \\ (1+q)x_1 & (1+q)x_2 \end{bmatrix}$, $q \neq 1$. The optimal value of the objective function $f(x_{optim}) = \begin{bmatrix} -0.3155 \\ -6.3092 \\ 1.8658 \end{bmatrix}$ for $x_{optim} = [0.0845\ 1.6908]^T$. The initial value is set to $x_0 = [0\ 0]^T$ and estimated results are given in Table 3.

Table 3: $x_1$, $x_2$, $f_1$, $f_2$ and $f_3$, for $q \in \{0.9, 0.95, 0.99, 1\}$

| $q$=1 Iterations = 14 norm=8.4562·10⁻⁷ | | | | | $q$=0.99 Iterations = 15 norm=3.7514·10⁻⁷ | | | | | $q$=0.95 Iterations =16 norm=4.8450·10⁻⁷ | | | | | $q$=0.9 Iterations = 17 norm=8.6438·10⁻⁷ | | | | |
|---|---|---|---|---|---|---|---|---|---|---|---|---|---|---|---|---|---|---|---|
| $x_1$ | $x_2$ | $f_1$ | $f_2$ | $f_3$ | $x_1$ | $x_2$ | $f_1$ | $f_2$ | $f_3$ | $x_1$ | $x_2$ | $f_1$ | $f_2$ | $f_3$ | $x_1$ | $x_2$ | $f_1$ | $f_2$ | $f_3$ |
| 0.4000 | 8.0000 | -0.4000 | -8.0000 | -1.0000 | 0.4000 | 8.0000 | -0.4000 | -8.0000 | -1.0000 | 0.4000 | 8.0000 | -0.4000 | -8.0000 | -1.0000 | 0.4000 | 8.0000 | -0.4000 | -8.0000 | -1.0000 |
| 0.2029 | 4.0581 | 0.0000 | 0.0000 | 63.1600 | 0.2029 | 4.0581 | 0.0000 | 0.0000 | 63.1600 | 0.1989 | 3.9779 | 0.0000 | 0.0000 | 63.1600 | 0.1937 | 3.8731 | 0.0000 | 0.0000 | 63.1600 |
| 0.1115 | 2.2306 | -0.1971 | -3.9419 | 15.5091 | 0.1115 | 2.2306 | -0.1971 | -3.9419 | 15.5091 | 0.1082 | 2.1633 | -0.2011 | -4.0221 | 14.8630 | 0.1040 | 2.0793 | -0.2063 | -4.1269 | 14.0387 |
| 0.0828 | 1.6556 | -0.2885 | -5.7694 | 3.9881 | 0.0828 | 1.6556 | -0.2885 | -5.7694 | 3.9881 | 0.0823 | 1.6465 | -0.2918 | -5.8367 | 3.6915 | 0.0822 | 1.6437 | -0.2960 | -5.9207 | 3.3345 |
| 0.0852 | 1.7049 | -0.3172 | -6.3444 | 1.7479 | 0.0852 | 1.7049 | -0.3172 | -6.3444 | 1.7479 | 0.0860 | 1.7205 | -0.3177 | -6.3535 | 1.7178 | 0.0869 | 1.7384 | -0.3178 | -6.3563 | 1.7084 |
| 0.0845 | 1.6891 | -0.3148 | -6.2951 | 1.9140 | 0.0845 | 1.6891 | -0.3148 | -6.2951 | 1.9140 | 0.0847 | 1.6943 | -0.3140 | -6.2795 | 1.9674 | 0.0851 | 1.7013 | -0.3131 | -6.2616 | 2.0295 |
| 0.0847 | 1.6938 | -0.3155 | -6.3109 | 1.8601 | 0.0847 | 1.6938 | -0.3155 | -6.3109 | 1.8601 | 0.0851 | 1.7028 | -0.3153 | -6.3057 | 1.8780 | 0.0857 | 1.7145 | -0.3149 | -6.2987 | 1.9018 |
| 0.0846 | 1.6924 | -0.3153 | -6.3062 | 1.8762 | 0.0846 | 1.6924 | -0.3153 | -6.3062 | 1.8762 | 0.0850 | 1.7000 | -0.3149 | -6.2972 | 1.9067 | 0.0855 | 1.7096 | -0.3143 | -6.2855 | 1.9467 |
| 0.0846 | 1.6928 | -0.3154 | -6.3076 | 1.8713 | 0.0846 | 1.6928 | -0.3154 | -6.3076 | 1.8713 | 0.0850 | 1.7009 | -0.3150 | -6.3000 | 1.8971 | 0.0856 | 1.7114 | -0.3145 | -6.2904 | 1.9301 |
| 0.0846 | 1.6927 | -0.3154 | -6.3072 | 1.8728 | 0.0846 | 1.6927 | -0.3154 | -6.3072 | 1.8728 | 0.0850 | 1.7006 | -0.3150 | -6.2991 | 1.9003 | 0.0855 | 1.7107 | -0.3144 | -6.2886 | 1.9362 |
| 0.0846 | 1.6927 | -0.3154 | -6.3073 | 1.8723 | 0.0846 | 1.6927 | -0.3154 | -6.3073 | 1.8723 | 0.0850 | 1.7007 | -0.3150 | -6.2994 | 1.8992 | 0.0855 | 1.7110 | -0.3145 | -6.2893 | 1.9340 |
| 0.0846 | 1.6927 | -0.3154 | -6.3073 | 1.8725 | 0.0846 | 1.6927 | -0.3154 | -6.3073 | 1.8725 | 0.0850 | 1.7007 | -0.3150 | -6.2993 | 1.8996 | 0.0855 | 1.7109 | -0.3145 | -6.2890 | 1.9348 |
| 0.0846 | 1.6927 | -0.3154 | -6.3073 | 1.8724 | 0.0846 | 1.6927 | -0.3154 | -6.3073 | 1.8724 | 0.0850 | 1.7007 | -0.3150 | -6.2993 | 1.8995 | 0.0855 | 1.7109 | -0.3145 | -6.2891 | 1.9345 |
| 0.0846 | 1.6927 | -0.3154 | -6.3073 | 1.8724 | 0.0846 | 1.6927 | -0.3154 | -6.3073 | 1.8724 | 0.0850 | 1.7007 | -0.3150 | -6.2993 | 1.8995 | 0.0855 | 1.7109 | -0.3145 | -6.2891 | 1.9346 |
| | | | | | 0.0846 | 1.6927 | -0.3154 | -6.3073 | 1.8724 | 0.0850 | 1.7007 | -0.3150 | -6.2993 | 1.8995 | 0.0855 | 1.7109 | -0.3145 | -6.2891 | 1.9345 |
| | | | | | | | | | | 0.0850 | 1.7007 | -0.3150 | -6.2993 | 1.8995 | 0.0855 | 1.7109 | -0.3145 | -6.2891 | 1.9346 |
| | | | | | | | | | | | | | | | 0.0855 | 1.7109 | -0.3145 | -6.2891 | 1.9346 |

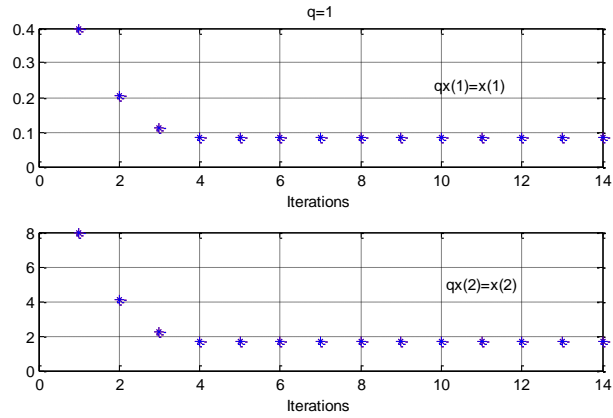

Figure 10: $x_1$ and $q*x_1$ (upwards) and $x_2$ and $q*x_2$ (downwards) (Iterations=14, $q$=1)

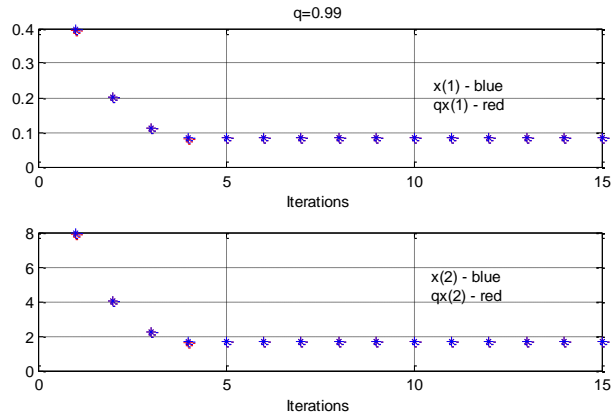

Figure 11: $x_1$ and $q*x_1$ (upwards) and $x_2$ and $q*x_2$ (downwards) (Iterations=15, $q$=0.99)

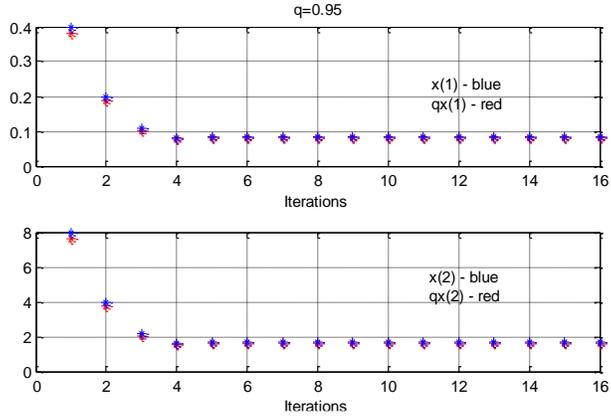

Figure 12: $x_1$ and $q*x_1$ (upwards) and $x_2$ and $q*x_2$ (downwards) (Iterations=16, $q$=0.95)

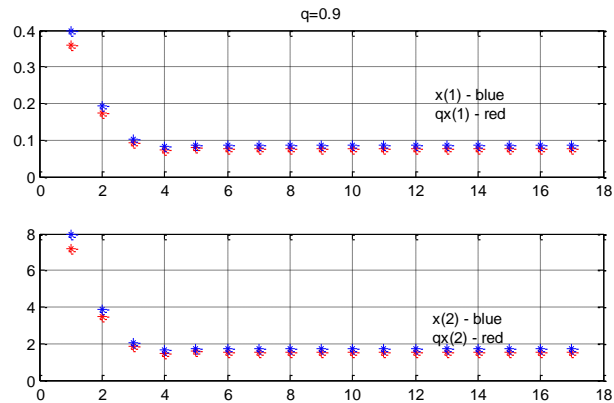

Figure 13: $x_1$ and $q*x_1$ (upwards) and $x_2$ and $q*x_2$ (downwards) (Iterations=17, $q$=0.9)

# 6 Conclusion

In this paper we proposed a novel Gauss-Newton method based on *q*-calculus, which solves nonlinear unconstrained optimization problems. The *q*-Gauss-Newton algorithm is the iterative procedure that converges only when the algorithm starts near to the optimal solution. Main advantage of the *q*-Gauss-Newton algorithm is that it does require only evaluations of the *q*-first-order derivatives of the objective function. The influence of the parameter *q* to the nonlinear problem solving is presented through three examples. The results show that the *q*-Gauss-Newton algorithm finds an optimal solution and speeds up the iterative procedure near the global minimum. In this way the *q*-Gauss-Newton method may be considered the basis of the *q*-damping algorithms designed without having to compute the *q*-Hessian matrix.